\documentclass[11pt]{article}
\usepackage{cite}
\usepackage{mathrsfs}
\usepackage{amssymb,amsfonts,amsmath}
\usepackage{enumerate}
\usepackage{indentfirst}
\usepackage{latexsym,bm}
\usepackage[noend]{algpseudocode}
\usepackage{algorithmicx,algorithm}
\usepackage{algorithm}
\usepackage{makeidx}
\usepackage{fancybox}
\usepackage{color}
\usepackage{multicol}
\usepackage[latin1]{inputenc}
\usepackage{graphicx}
\usepackage{epstopdf}
\usepackage{pstricks,pst-node,pst-text,pst-3d}
\usepackage{tikz}
\newtheorem{thm}{Theorem}[section]
\newtheorem{cor}[thm]{Corollary}
\newtheorem{lem}[thm]{Lemma}
\newtheorem{op}[thm]{Problem}

\newtheorem{pro}[thm]{Proposition}

\newenvironment{pf}{{\noindent \it \bf Proof:}}{{\hfill$\Box$}\\}

\newcommand{\2}{\vspace{0.2cm}}

\begin{document}

\title{\bf Directed Steiner tree packing and directed tree connectivity}
\author{Yuefang Sun$^{1}$\thanks{E-mail: yuefangsun2013@163.com. Research supported by Zhejiang Provincial Natural Science Foundation (No. LY20A010013).} \mbox{ } and  \mbox{ }
Anders Yeo${}^{2,3}$\thanks{E-mail: yeo@imada.sdu.dk. Research supported by the Independent
Research Fund Denmark under grant number DFF 7014-00037B} \\ \\
$^{1}$ School of Mathematics and Statistics, Ningbo University,\\
Zhejiang 315211, P. R. China\\ \\
$^{2}$ Department of Mathematics and Computer Science, \\
University of Southern Denmark, Odense, Denmark \\ \\
$^{3}$ Department of Mathematics and Applied Mathematics\\
University of Johannesburg, Auckland Park, 2006 South Africa }
\date{}
\maketitle

\begin{abstract}
For a digraph $D=(V(D), A(D))$, and a set $S\subseteq V(D)$ with
$r\in S$ and $|S|\geq 2$, an $(S, r)$-tree is an out-tree $T$ rooted at $r$ with
$S\subseteq V(T)$. Two $(S, r)$-trees
$T_1$ and $T_2$ are said to be arc-disjoint if $A(T_1)\cap
A(T_2)=\emptyset$. Two arc-disjoint $(S, r)$-trees $T_1$ and $T_2$ are
said to be internally disjoint if $V(T_1)\cap V(T_2)=S$. Let $\kappa_{S,r}(D)$ and $\lambda_{S,r}(D)$ be the maximum
number of internally disjoint and arc-disjoint $(S, r)$-trees in $D$,
respectively. The generalized $k$-vertex-strong connectivity
of $D$ is defined as
$$\kappa_k(D)= \min \{\kappa_{S,r}(D)\mid S\subset V(D), |S|=k, r\in
S\}.$$ Similarly, the generalized $k$-arc-strong connectivity
of $D$ is defined as
$$\lambda_k(D)= \min \{\lambda_{S,r}(D)\mid S\subset V(D), |S|=k, r\in
S\}.$$
The generalized $k$-vertex-strong connectivity and generalized $k$-arc-strong
connectivity are also called directed tree connectivity which extends the well-established tree connectivity on undirected graphs to directed graphs and could be seen as a generalization of classical connectivity of digraphs.

In this paper, we completely determine the complexity for both $\kappa_{S, r}(D)$ and $\lambda_{S, r}(D)$ on general digraphs, symmetric digraphs and Eulerian digraphs. In particular, among our results, we prove and use the NP-completeness of 2-linkage problem restricted to Eulerian digraphs. We also give sharp bounds and characterizations for the two parameters $\kappa_k(D)$ and $\lambda_k(D)$.
\vspace{0.3cm}
\vspace{0.3cm}
\end{abstract}


\section{Introduction}\label{sec:intro}

We refer the readers to \cite{Bang-Jensen-Gutin, Bang-Jensen-Gutin2,
Bondy} for graph theoretical notation and terminology not given
here\footnote{Note that all digraphs considered in this paper have no parallel arcs or loops.}.
For a graph $G=(V,E)$ and a set $S\subseteq V$ of at
least two vertices, an {\em $S$-Steiner tree} or, simply, an {\em
$S$-tree} is a tree $T$ of $G$ with $S\subseteq V(T)$. Two $S$-trees $T_1$ and $T_2$ are said to be {\em
edge-disjoint} if $E(T_1)\cap E(T_2)=\emptyset$. Two edge-disjoint
$S$-trees $T_1$ and $T_2$ are said to be {\em internally disjoint}
if $V(T_1)\cap V(T_2)=S$. The basic problem of {\sc Steiner Tree
Packing} is defined as follows: the input consists of an undirected
graph $G$ and a subset of vertices $S\subseteq V(D)$, the goal is to
find a largest collection of edge-disjoint $S$-Steiner trees.  The Steiner tree packing problem has applications in VLSI circuit design \cite{Grotschel-Martin-Weismantel, Sherwani}. 
In this application, a Steiner tree is needed to share an electronic signal by a set of terminal nodes. Another application,
arises in the Internet Domain \cite{Li-Mao5}. Imagine that a given graph $G$
represents a network. We choose arbitrary $k$ vertices as nodes. Suppose one of them is a {\em broadcaster}, and all other nodes are either {\em users} or {\em routers} (also called {\em switches}). The broadcaster wants to broadcast as many streams of movies as possible, so that the users have
the maximum number of choices. Each stream of movie is broadcasted via a tree connecting
all the users and the broadcaster. In essence we need to find the maximum number Steiner
trees connecting all the users and the broadcaster. Clearly, it is a Steiner tree packing problem.
Besides the classical version, people also study some other
variations, such as packing internally disjoint Steiner trees,
packing directed Steiner trees and packing strong subgraphs \cite{Cheriyan-Salavatipour,
DeVos-McDonald-Pivotto, Kriesell, Lau, Sun-Gutin-Yeo-Zhang, West-Wu}.

The generalized $k$-connectivity $\kappa_k(G)$ of a graph $G=(V,E)$
which is related to internally disjoint Steiner trees packing
problem was introduced by Hager \cite{Hager} in 1985 ($2\le k\le
|V|$). The {\em generalized local connectivity} $\kappa_S(G)$ is the
maximum number of internally disjoint $S$-trees in $G$, and the {\em
generalized $k$-connectivity} of $G$ is defined as
$$\kappa_k(G)=\min\{\kappa_S(G)\mid S\subseteq V(G), |S|=k\}.$$
Li, Mao and Sun \cite{Li-Mao-Sun} introduced the following concept
of generalized $k$-edge-connectivity which is related to basic
problem of Steiner tree packing. The {\em generalized local
edge-connectivity} $\lambda_S(G)$ is the maximum number of
edge-disjoint $S$-trees in $G$ and the {\em generalized
$k$-edge-connectivity} is defined as
$$\lambda_k(G)=\min\{\lambda_S(G)\mid S\subseteq V(G), |S|=k\}.$$
Observe that $\kappa_2(G)=\kappa(G)$ and $\lambda_2(G)=\lambda(G)$.
Hence, these two parameters could be seen as natural generalizations
of connectivity and edge connectivity of a graph, respectively.
This type of generalized connectivity is also called {\em tree connectivity} in
the literature, and has become an established area in graph theory,
see a recent monograph \cite{Li-Mao5} by Li and Mao on this topic.

Below we summarize the known complexity for determining whether $\kappa_S(G) \geq \ell$
and for determining whether $\lambda_S(G) \geq \ell$ for any $k$ and $\ell$.
Note that if $k=2$ or if $\ell=1$, then the problem is equivalent to normal connectivity and therefore polynomial. In all
other cases the following holds.

\begin{center}
\begin{tabular}{|c||c|c|c|} \hline
$\lambda_S(G) \geq \ell$? & $k=3$                     & $k \geq 4$             &  $k$ part \\
$ |S|=k$                 &                            & constant               &  of input \\ \hline  \hline
$\ell =2$                & Polynomial \cite{Chen-Li-Liu-Mao}       & Polynomial \cite{Chen-Li-Liu-Mao}   &  NP-complete \cite{Chen-Li-Liu-Mao} \\ \hline
$\ell \geq 3$ constant   & Polynomial \cite{Chen-Li-Liu-Mao}       & Polynomial \cite{Chen-Li-Liu-Mao}   &  NP-complete \cite{Chen-Li-Liu-Mao}  \\ \hline
$\ell$ part of input     & NP-complete \cite{Chen-Li-Liu-Mao}      & NP-complete \cite{Chen-Li-Liu-Mao}  &  NP-complete \cite{Chen-Li-Liu-Mao} \\ \hline
\end{tabular}
\end{center}

\begin{center}
\begin{tabular}{|c||c|c|c|} \hline
$\kappa_S(G) \geq \ell$? & $k=3$                              & $k \geq 4$                &  $k$ part \\
$ |S|=k$                 &                                    & constant                  &  of input \\ \hline  \hline
$\ell =2$                & Polynomial \cite{Li-Li-Zhou}       & Polynomial \cite{Li-Li}   &  NP-complete \cite{Li-Li} \\ \hline
$\ell \geq 3$ constant   & Polynomial \cite{Li-Li-Zhou}       & Polynomial \cite{Li-Li}   &  NP-complete \cite{Li-Li}  \\ \hline
$\ell$ part of input     & NP-complete \cite{Chen-Li-Liu-Mao} & NP-complete \cite{Li-Li}  &  NP-complete \cite{Li-Li} \\ \hline
\end{tabular}
\end{center}

We will now prove similar results for directed graphs.
An {\em out-tree (respectively, in-tree)} is an oriented tree in which every vertex
except one, called the {\em root}, has in-degree (respectively, out-degree) one.
An {\em out-branching (respectively, in-branching)} of $D$ is a spanning out-tree
(respectively, in-tree) in $D$. For a digraph $D=(V(D), A(D))$, and a set $S\subseteq V(D)$ with
$r\in S$ and $|S|\geq 2$, a {\em directed $(S, r)$-Steiner tree} or,
simply, an {\em $(S, r)$-tree} is an out-tree $T$ rooted at $r$ with
$S\subseteq V(T)$ \cite{Cheriyan-Salavatipour}. Two $(S, r)$-trees
$T_1$ and $T_2$ are said to be {\em arc-disjoint} if $A(T_1)\cap
A(T_2)=\emptyset$. Two arc-disjoint $(S, r)$-trees $T_1$ and $T_2$ are
said to be {\em internally disjoint} if $V(T_1)\cap V(T_2)=S$.

Cheriyan and Salavatipour \cite{Cheriyan-Salavatipour} introduced
and studied the following two directed Steiner tree packing
problems. {\sc Arc-disjoint directed Steiner tree packing} (ADSTP):
The input consists of a digraph $D$ and a subset of vertices
$S\subseteq V(D)$ with a root $r$, the goal is to find a largest
collection of arc-disjoint $(S, r)$-trees. {\sc Internally-disjoint
directed Steiner tree packing} (IDSTP): The input consists of a
digraph $D$ and a subset of vertices $S\subseteq V(D)$ with a root
$r$, the goal is to find a largest collection of internally disjoint
$(S, r)$-trees.

The following concept of directed tree connectivity is
related to directed Steiner tree packing problem and is a natural
extension of tree connectivity of undirected graphs to directed
graphs. Let $\kappa_{S,r}(D)$ (respectively, $\lambda_{S,r}(D)$) be the maximum
number of internally disjoint (respectively, arc-disjoint) $(S, r)$-trees in $D$. The {\em generalized $k$-vertex-strong connectivity}
of $D$ is defined as
$$\kappa_k(D)= \min \{\kappa_{S,r}(D)\mid S\subset V(D), |S|=k, r\in
S\}.$$ Similarly, the {\em generalized $k$-arc-strong connectivity}
of $D$ is defined as
$$\lambda_k(D)= \min \{\lambda_{S,r}(D)\mid S\subset V(D), |S|=k, r\in
S\}.$$ By definition, when $k=2$, $\kappa_2(D)=\kappa(D)$ and
$\lambda_2(D)=\lambda(D)$. Hence, these two parameters could be seen
as generalizations of vertex-strong connectivity and arc-strong
connectivity of a digraph. The generalized $k$-vertex-strong connectivity and $k$-arc-strong
connectivity are also called {\em directed tree connectivity}.

We will in this paper prove the non-cited results of the following tables.

\begin{center}
\begin{tabular}{|c||c|c|c|} \hline
\multicolumn{4}{|c|}{Table 1: Directed graphs} \\ \hline
$\lambda_{S,r}(D) \geq \ell$? & $k=3$                                    & $k \geq 4$     &  $k$ part \\
$ |S|=k$                     &                                          & constant       &  of input \\ \hline  \hline
$\ell =2$                    & NP-complete \cite{Cheriyan-Salavatipour} & NP-complete    &  NP-complete  \\ \hline
$\ell \geq 3$ constant       & NP-complete                              & NP-complete    &  NP-complete  \\ \hline
$\ell$ part of input         & NP-complete                              & NP-complete    &  NP-complete  \\ \hline
\end{tabular}
\end{center}

\begin{center}
\begin{tabular}{|c||c|c|c|} \hline
\multicolumn{4}{|c|}{Table 2: Directed graphs} \\ \hline
$\kappa_{S,r}(D) \geq \ell$? & $k=3$                                    & $k \geq 4$     &  $k$ part \\
$ |S|=k$                     &                                          & constant       &  of input \\ \hline  \hline
$\ell =2$                    & NP-complete \cite{Cheriyan-Salavatipour} & NP-complete    &  NP-complete  \\ \hline
$\ell \geq 3$ constant       & NP-complete                              & NP-complete    &  NP-complete  \\ \hline
$\ell$ part of input         & NP-complete                              & NP-complete    &  NP-complete  \\ \hline
\end{tabular}
\end{center}

A digraph $D$ is {\em symmetric} if every arc in $D$ belongs to a $2$-cycle.  That is, if $xy \in A(D)$ then $yx \in A(D)$.
 In other words, a symmetric
digraph $D$ can be obtained from its underlying undirected graph $G$
by replacing each edge of $G$ with the corresponding arcs of both
directions, that is, $D=\overleftrightarrow{G}.$
For a digraph $D$, its {\em reverse} $D^{\rm rev}$ is a digraph with the same vertex set such that $xy\in A(D^{\rm rev})$ if and only if $yx\in A(D)$.
Note that if a digraph $D$ is {\em symmetric} then $D^{\rm rev}=D$.
We will in this paper also prove the results in the following tables.

\begin{center}
\begin{tabular}{|c||c|c|c|} \hline
\multicolumn{4}{|c|}{Table 3: Symmetric digraphs} \\ \hline
$\lambda_{S,r}(D) \geq \ell$? & $k=3$        & $k \geq 4$     &  $k$ part \\
$ |S|=k$                      &              & constant       &  of input \\ \hline  \hline
$\ell =2$                     & Polynomial   & Polynomial     &  Polynomial   \\ \hline
$\ell \geq 3$ constant        & Polynomial   & Polynomial     &  Polynomial   \\ \hline
$\ell$ part of input          & Polynomial   & Polynomial     &  Polynomial   \\ \hline
\end{tabular}
\end{center}

\begin{center}
\begin{tabular}{|c||c|c|c|} \hline
\multicolumn{4}{|c|}{Table 4: Symmetric digraphs} \\ \hline
$\kappa_{S,r}(D) \geq \ell$? & $k=3$        & $k \geq 4$     &  $k$ part \\
$ |S|=k$                     &              & constant       &  of input \\ \hline  \hline
$\ell =2$                    & Polynomial   & Polynomial     &  NP-complete  \\ \hline
$\ell \geq 3$ constant       & Polynomial   & Polynomial     &  NP-complete  \\ \hline
$\ell$ part of input         & NP-complete  & NP-complete    &  NP-complete  \\ \hline
\end{tabular}
\end{center}

In a digraph $D$, the maximum number of arc-disjoint $(x,y)$-paths is denoted by $\lambda_D(x,y)$, or just by $\lambda(x,y)$ if $D$ is clear from the context.
In order to determine the complexities of Table~3 we will actually prove the following stronger result.

\begin{thm}\label{LambdaEulerian}
If $D$ is an Eulerian digraph and $S \subseteq V(D)$ and $r \in S$, then $\lambda_{S,r}(D) \geq \ell$ if and only if $\lambda_D(r,s) \geq \ell$ for all $s \in S \setminus \{r\}$.
\end{thm}

As one can determine $\lambda_D(r,s)$ in polynomial time for any $r$ and $s$ in $D$ we note that
Theorem~\ref{LambdaEulerian} implies that we can extend the results in Table~3 from symmetric digraphs to Eulerian digraphs implying that the following also holds.

\begin{center}
\begin{tabular}{|c||c|c|c|} \hline
\multicolumn{4}{|c|}{Table 5: Eulerian digraphs} \\ \hline
$\lambda_{S,r}(D) \geq \ell$? & $k=3$        & $k \geq 4$     &  $k$ part \\
$ |S|=k$                      &              & constant       &  of input \\ \hline  \hline
$\ell =2$                     & Polynomial   & Polynomial     &  Polynomial   \\ \hline
$\ell \geq 3$ constant        & Polynomial   & Polynomial     &  Polynomial   \\ \hline
$\ell$ part of input          & Polynomial   & Polynomial     &  Polynomial   \\ \hline
\end{tabular}
\end{center}

If we consider $\kappa_{S,r}$ instead of $\lambda_{S,r}$ for Eulerian digraphs then
we will prove the following results.

\begin{center}
\begin{tabular}{|c||c|c|c|} \hline
\multicolumn{4}{|c|}{Table 6: Eulerian digraphs} \\ \hline
$\kappa_{S,r}(D) \geq \ell$? & $k=3$        & $k \geq 4$     &  $k$ part \\
$ |S|=k$                      &              & constant       &  of input \\ \hline  \hline
$\ell =2$                     & NP-complete             & NP-complete               &  NP-complete   \\ \hline
$\ell \geq 3$ constant        & NP-complete             & NP-complete               &  NP-complete   \\ \hline
$\ell$ part of input          & NP-complete  & NP-complete    &  NP-complete   \\ \hline
\end{tabular}
\end{center}

It may be slightly surprising that for Eulerian digraphs the complexity of deciding if $\lambda_{S,r}(D) \geq \ell$ is always polynomial,
while the complexity of deciding if $\kappa_{S,r}(D) \geq \ell$ is always NP-complete.


In this paper we furthermore prove some inequalities  for directed tree connectivity. Let $D$
be a strong digraph of order $n$. For $2\leq k\leq n$, we prove
that $1\leq \kappa_k(D)\leq n-1$ and $1\leq \lambda_k(D)\leq n-1$
(Theorem~\ref{thma}). All bounds are sharp, we also characterize
those digraphs $D$ for which $\kappa_k(D)$~(respectively, $\lambda_k(D)$) attains
the upper bound.       
We also study the relation between the
directed tree connectivity and classical connectivity of digraphs by
showing that $\kappa_k(D)\leq \kappa(D)$ and $\lambda_k(D)\leq
\lambda(D)$ (Theorem~\ref{thmb}). Furthermore, these bounds are sharp. In
Section~\ref{sec:bounds}, the sharp Nordhaus-Gaddum type bounds for
$\lambda_k(D)$ are also given; moreover, we characterize those
extremal digraphs for the lower bounds (Theorem~\ref{thmf}). 


\paragraph{Additional Terminology and Notation.}
A digraph $D$ is {\em semicomplete} if for
every distinct $x,y\in V(D)$ at least one of the arcs $xy,yx$ is in
$D$. For a digraph $D=(V(D), A(D))$, the {\em complement digraph},
denoted by $D^c$, is a digraph with vertex set $V(D^c)=V(D)$ such
that $xy\in A(D^c)$ if and only if $xy\not\in A(D)$. 
\section{Complexity for packing directed Steiner trees}\label{sec:complexity}

Let $D$ be a digraph and $S \subseteq V(D)$ with $|S|=k$. It is
natural to consider the following problem: what is the complexity of
deciding whether $\kappa_{S, r}(D)\ge \ell$~(respectively, $\lambda_{S, r}(D)\ge
\ell)$? where $r\in S$ is a root. If $k=2$, say $S=\{r, x\}$, then
the problem of deciding whether $\kappa_{S, r}(D)\ge
\ell$~(respectively, $\lambda_{S, r}(D)\ge \ell)$ is equivalent to deciding whether
$\kappa(r, x)\ge \ell$~(respectively, $\lambda(r, x)\ge \ell)$, and so is
polynomial-time solvable (see \cite{Bang-Jensen-Gutin}), where
$\kappa(r, x)$~(respectively, $\lambda(r, x))$ is the local vertex-strong
(respectively, arc-strong) connectivity from $r$ to $x$. If $\ell=1$, then the above problem is also polynomial-time solvable by the well-known fact that every strong digraph has an out- and in-branching rooted at any vertex, and these branchings can be found in polynomial-time. Hence, it remains to consider the case that $k\geq 3,
\ell\geq 2$.
That is, we will prove the results in Tables~1 and 2. The known parts of Table~1 and 2 are the following.

\begin{thm}\label{thm04}\cite{Cheriyan-Salavatipour}
Let $D$ be a digraph and $S \subseteq V(D)$ with $|S|=3$. The
problem of deciding whether $\kappa_{S, r}(D)\ge 2$ is NP-hard,
where $r\in S$.
\end{thm}

\begin{thm}\label{thm05}\cite{Cheriyan-Salavatipour}
Let $D$ be a digraph and $S \subseteq V(D)$ with $|S|=3$. The
problem of deciding whether $\lambda_{S, r}(D)\ge 2$ is NP-hard,
where $r\in S$.
\end{thm}

We will now extend the above results to the following.

\begin{thm}\label{thmd}
Let $k\geq 3$ and $\ell\geq 2$ be fixed integers (considered as constants). Let $D$ be a
digraph and $S \subseteq V(D)$ with $|S|=k$ and $r \in S$. Both the following problems are NP-complete.
\begin{itemize}
\item Is $\kappa_{S, r}(D)\ge \ell$?
\item Is $\lambda_{S, r}(D)\ge \ell$?
\end{itemize}
\end{thm}

\begin{figure}
\begin{center}
\tikzstyle{vertexB}=[circle,draw, minimum size=15pt, scale=0.7, inner sep=0.1pt]
\tikzstyle{vertexL}=[circle,draw, minimum size=15pt, scale=0.9, inner sep=0.1pt]
\begin{tikzpicture}[scale=0.25]
\node (r) at (4,6) [vertexL] {$r$};
\node (s1) at (32,4) [vertexL] {$s_1$};
\node (s2) at (32,8) [vertexL] {$s_2$};
\draw (1,1) rectangle (36,12);
\draw (30,2) rectangle (34,10);

\draw (30,15) rectangle (34,28);
\draw (8,15) rectangle (12,28);
\draw (17,15) rectangle (21,28);

\node (s3) at (32,17) [vertexL] {$s_3$};
\node (s4) at (32,20) [vertexL] {$s_4$};
\node (sk) at (32,26) [vertexB] {$s_{k-1}$};

\node (u1) at (10,17) [vertexL] {$u_1$};
\node (w1) at (19,17) [vertexL] {$w_1$};

\node (u2) at (10,20) [vertexL] {$u_2$};
\node (w2) at (19,20) [vertexL] {$w_2$};

\node (ul) at (10,26) [vertexB] {$u_{\ell-2}$};
\node (wl) at (19,26) [vertexB] {$w_{\ell-2}$};

\node () at (10,30) [scale=1.1] {$U$};
\node () at (19,30) [scale=1.1] {$W$};

\node () at (18,6) [scale=1.5] {$D^*$};
\node () at (32,23) [scale=0.8] {$\vdots$};
\node () at (10,23) [scale=0.8] {$\vdots$};
\node () at (19,23) [scale=0.8] {$\vdots$};

\draw [->, line width=0.02cm] (r) to (8,15);
\draw [->, line width=0.02cm] (u1) to (w1);
\draw [->, line width=0.02cm] (u2) to (w2);
\draw [->, line width=0.02cm] (ul) to (wl);
\draw [->, line width=0.02cm] (21,22) to (30,22);
\draw [->, line width=0.02cm] (21,15) to (30,10);
\draw [->, line width=0.02cm] (32,10) to (32,15);

\end{tikzpicture}
\end{center}
\caption{In the above digraph, $D$, the following statements hold. \hspace{3cm}
(1): $\kappa_{\{r,s_1,s_2,\ldots,s_{k-1}\}, r}(D)\geq \ell$  if and only if $\kappa_{\{r,s_1,s_2\},r}(D^*) \geq 2$. \mbox{ } \hspace{0.7cm} \mbox{ } \hspace{0.1cm}
(2): $\lambda_{\{r,s_1,s_2,\ldots,s_{k-1}\}, r}(D)\geq \ell$ if and only if $\lambda_{\{r,s_1,s_2\},r}(D^*) \geq 2$.}  \label{pic1}
\end{figure}

\begin{pf}
We first consider the problem of determining if $\kappa_{S, r}(D)\ge \ell$.
It is easy to see that this problem belongs to NP. To show it is NP-hard,
we reduce from the case when $k=3$ and $\ell=2$.
That is let $D^*$ be a digraph, $S^* \subseteq V(D^*)$ with $|S^*|=3$ and $r \in S^*$ where
we want to determine if $\kappa_{S^*, r}(D)\ge 2$.
This problem is NP-hard by Theorem~\ref{thm04}.

We will construct a new digraph $D$ containing a set $S \subseteq V(D)$ with $|S|=k$ and $r \in S$,
such that $\kappa_{S, r}(D)\geq \ell$ if and only if $\kappa_{S^*,r}(D^*) \geq 2$. This would complete the proof
for the  problem of determining if $\kappa_{S, r}(D)\ge \ell$.
Assume that $S^*=\{r,s_1,s_2\}$ and let $V(D)=V(D^*) \cup U \cup  W \cup \{s_3,s_4,\ldots,s_{k-1}\}$, where
$U=\{u_1,u_2,\ldots,u_{\ell-2}\}$ and $W=\{w_1,w_2,\ldots,w_{\ell-2}\}$.

Let $S=\{r,s_1,s_2,s_3,\ldots,s_{k-1}\}$ and
let the arc-set of $D$ be the following (see Figure~\ref{pic1}).

\[
\begin{array}{rclll}
A(D) & = & A(D^*) & \cup & \{ru_i, u_i w_i, w_i s \; | \; i=1,2,\ldots , \ell-2 \mbox{ and } s \in S-r \} \\
 & & & \cup & \{s_i s_j \; | \; i=1,2 \mbox{ and } j=3,4,\ldots,k-1 \} \\
\end{array}
\]

First assume that  $\kappa_{S^*,r}(D^*) \geq 2$ and let $T_1^*$ and $T_2^*$ be the two internally disjoint $(S^*,r)$-trees in $D^*$.
Add all arcs from $s_1$ to $\{s_3,s_4,\ldots,s_{k-1}\}$ to $T_1^*$ and add all arcs from
$s_2$ to $\{s_3,s_4,\ldots,s_{k-1}\}$ to $T_2^*$ in order to obtain two internally disjoint $(S,r)$-trees in $D$.
Furthermore, for $i=1,2,\ldots,\ell-2$, we note that the path $r u_i w_i$ together with all arcs from $w_i$ to $S-r$ gives us an
$(S,r)$-tree in $D$. It is not difficult to see that the $\ell$ constructed $(S,r)$-trees are all internally disjoint
and therefore $\kappa_{S, r}(D)\ge \ell$.

Conversely, if $\kappa_{S, r}(D)\ge \ell$, then let $T_1,T_2,\ldots,T_{\ell}$ be the $\ell$ internally disjoint $(S,r)$-trees in $D$.
If we remove all the $T_i$'s containing a vertex from $U$ then we are left with at least two trees, not containing any vertex from $U$ or
$W$ (as after removing $U$ the vertices in $W$ have no arc into them). Removing all vertices from $\{s_3,s_4,\ldots,s_{k-2}\}$ from
the two remaining $T_i$'s gives us two internally disjoint $(S^*,r)$-trees in $D^*$ and therefore $\kappa_{S^*,r}(D^*) \geq 2$.

This completes the proof for determining if $\kappa_{S, r}(D)\ge \ell$.
The case when we want to determine if $\lambda_{S, r}(D)\ge \ell$ is similar.
We use exactly the same reduction, except in $D^*$ we look at the problem of determining if $\lambda_{S^*, r}(D^*)\ge 2$, which
is also NP-hard by Theorem~\ref{thm05}. Otherwise the reduction is identical (but delete the arcs $u_iw_i$ instead of the vertices in $U$ in the
above proof).
\end{pf}

The above results imply the entries in Tables~1 and 2.

\section{Symmetric digraphs}

Now we turn our attention to symmetric digraphs. We first need the following theorem by Robertson and Seymour.

\begin{thm}\label{LinkageUndirected} \cite{Robertson-Seymour}
Let $G$ be a graph and let $s_1,s_2,\ldots,s_r, t_1,t_2,\ldots,t_r$ be $2r$ disjoint vertices in $G$.
We can in $O(|V(G)|^3)$ time decide if there exists an $(s_i,t_i)$-path, $P_i$, such that all $P_1,P_2,\ldots,P_r$ are vertex disjoint.
\end{thm}

\begin{cor}\label{corLinkageUndirected}
Let $D$ be a symmetric digraph and let $s_1,s_2,\ldots,s_r, t_1,t_2,\ldots,t_r$ be vertices in $D$ (not necessarily disjoint) and
let $S \subseteq V(D)$. We can in $O(|V(G)|^3)$ time decide if there for all $i=1,2,\ldots,r$ exists an $(s_i,t_i)$-path, $P_i$,
such that no internal vertex of any $P_i$ belongs to $S$ or to any path $P_j$ with $j \not=i$ (the end-points of $P_j$ can also not
be internal vertices of $P_i$).
\end{cor}
\begin{pf}
Let $D$ and $S$ and  $s_1,s_2,\ldots,s_r, t_1,t_2,\ldots,t_r$ be defined as in the corollary.
If some vertex, $x$, appears $q$ times in the sequence $s_1,s_2,\ldots,s_r, t_1,t_2,\ldots,t_r$, where $q \geq 1$, then
we make $q$ copies, $x_1,x_2, \ldots, x_q$ of $x$ and replace $x$ with these $q$ copies.  All arcs in $D$ that entered $x$ now enter
every copy of $x$ and all arcs out of $x$ in $D$ now go out of every copy of $x$.
We then replace the $q$ copies of $x$ in the sequence by a separate copy of $x$.
Finally, every vertex in $S$, that does not appear in the sequence $s_1,s_2,\ldots,s_r, t_1,t_2,\ldots,t_r$ gets deleted.
This results in a new symmetric digraph $D'$ and a sequence $s_1',s_2',\ldots,s_r', t_1',t_2',\ldots,t_r'$, of $2r$ disjoint vertices in $D'$.
Note that every vertex in $S$ only appears as vertices in the sequence $s_1',s_2',\ldots,s_r', t_1',t_2',\ldots,t_r'$ (or gets
removed completely if it is not in this sequence).  Let $G'$ be the underlying graph of $D'$.

If $P'$ is an $(s_i',t_i')$-path in $G'$, then $P'$ corresponds to an $(s_i',t_i')$-path in $D'$ and
therefore also to an $(s_i,t_i)$-path in $D$.
If $P_1',P_2',\ldots,P_r'$ are vertex-disjoint paths in $G'$, where $P_i'$ is an $(s_i',t_i')$-path in $G'$,
then the corresponding paths in $D$ are the desired paths we were looking for.

Conversely, if for all $i=1,2,\ldots,r$ there exists an $(s_i,t_i)$-path, $P_i$, in $D$,
such that no internal vertex of any $P_i$ belongs to $S$ or to any path $P_j$ with $j \not=i$,
then the corresponding paths in $G'$ are all vertex-disjoint
$(s_i',t_i')$-paths. We are now done by Theorem~\ref{LinkageUndirected}.
\end{pf}

\begin{thm}\label{thmPOLsym}
Let $k\geq 3$ and $\ell \geq 2$ be fixed integers and let $D$ be a symmetric digraph. Let $S \subseteq V(D)$ with $|S|=k$ and let $r$ be an arbitrary vertex in $S$. Let $A_0,A_1,A_2,\ldots , A_{\ell}$ be a partition of the arcs in $D[S]$.

We can in time $O(n^{\ell k - 2 \ell + 3} \cdot (2k-3)^{\ell (2k-3)})$
decide if there exist $\ell$ internally disjoint $(S,r)$-trees, $T_1,T_2,\ldots ,T_{\ell}$, with $A(T_i) \cap A[S] = A_i$ for all $i=1,2,\ldots,\ell$ (note that $A_0$ are the arcs in $D[S]$ not used in any of the trees).
\end{thm}
\begin{pf}
Let $T$ be any $(S,r)$-tree in $D$. Let $R$ be all vertices in $T$ with degree (in the underlying graph of $T$) at least three and let $L$
be all the leaves of $T$ (in the underlying graph of $T$). It is well-known that $|R| \leq |L|-2$. Furthermore we will assume that all
$(S,r)$-trees  considered are minimal (i.e if we delete a vertex in the tree it is not an $(S,r)$-tree anymore), which implies that all
leaves are in $S$.  Under this assumption we have $|R| \leq |L|-2 \leq |S|-2$. Note that every vertex $x \in V(T) \setminus (R \cup S)$ has
$d_T^+(x)=d_T^-(x)=1$. We now define the {\em skeleton} of $T$ as the tree we obtain from $T$ by by-passing all vertices in $ V(T) \setminus (R \cup S)$
(i.e. if $ux,xv \in A(T)$ and  $x \in V(T) \setminus (R \cup S)$ then replace $ux$ and $xv$ by $uv$).

Let $T^s$ be a skeleton of an $(S,r)$-tree in $D$ and recall that $V(T^s)$ consists of $S$ as well as at most $|S|-2$ ($=k-2$) other vertices.
Therefore there are less than $n^{k-2}$ possibilities for $V(T^s)$.
Once $V(T^s)$ is known there are at most $(2k-3)^{2k-3}$ possibilities for the arcs of $T^s$ (as for each $y \in V(T^s)\setminus \{r\}$ we need
to pick the vertex with an arc into $y$ and there are at most $2k-3$ possibilities for picking this vertex).
This implies that there are at most $n^{k-2} \cdot (2k-3)^{2k-3}$ different skeletons of minimal $(S,r)$-trees in $D$.

Our algorithm will try all possible $\ell$-tuples, ${\cal T}^s = (T_1^s, T_2^s,\ldots,T_{\ell}^s)$, of skeletons of $(S,r)$-trees and determine
if there are $\ell$ internally disjoint $(S,r)$-trees, $T_1,T_2,\ldots ,T_{\ell}$,
with $A(T_i) \cap A[S] = A_i$ for all $i=1,2,\ldots,\ell$ and such that $T_i^s$ is the skeleton of $T_i$.
If such a set of trees exists for any ${\cal T}^s$, then we return this solution and if no such set of trees exist for any ${\cal T}^s$ we return
that no solution exists.
We will first prove that this algorithm gives the correct answer and then compute its time complexity.

If our algorithm returns a solution, then clearly a solution exists. So now assume that a solution exists and
let $T_1,T_2,\ldots ,T_{\ell}$ be the desired set of internally disjoint $(S,r)$-trees.
When we consider ${\cal T}^s = (T_1^s, T_2^s,\ldots,T_{\ell}^s)$, where $T_i^s$ is the skeleton of $T_i$, our algorithm will find a
solution, so the algorithm always returns a solution if one exists.

We will now analyse the time complexity. The number of different $\ell$-tuples, ${\cal T}^s$, that we need to consider is bounded by the following.

\[
\left( n^{k-2} \cdot (2k-3)^{2k-3} \right)^{\ell} = n^{\ell (k-2)} \cdot (2k-3)^{\ell (2k-3)}
\]

Given such a $\ell$-tuples, ${\cal T}^s$, we need to determine if there exist $\ell$ internally disjoint $(S,r)$-trees, $T_1,T_2,\ldots ,T_{\ell}$,
with $A(T_i) \cap A[S] = A_i$ for all $i=1,2,\ldots,\ell$ and such that $T_i^s$ is the skeleton of $T_i$.
We first check that the arcs in $A_i$ belong to the skeleton $T_i^s$ and that no vertex in $V(D)\setminus S$ belongs to more than one skeleton.
If the above does not hold then the desired trees do not exist, so assume the above holds. We will now use Corollary~\ref{corLinkageUndirected}.
In fact, for every arc $uv \not\in A(D[S])$ that belongs to some skeleton $T_i^s$ we want to find a $(u,v)$-path in $D - A(D[S])$,
such that no internal vertex on any path belongs to $S$ or to a different path.
By Corollary~\ref{corLinkageUndirected} this can be done in $O(n^3)$ time.  If such paths exist, then substituting the arcs $uv$ by the $(u,v)$-paths
we obtain the desired  $(S,r)$-trees. And if the paths do not exist the desired  $(S,r)$-trees do not exist.
Therefore the algorithm works correctly and has complexity $O(n^{\ell k-2 \ell + 3} \cdot (2k-3)^{\ell (2k-3)} )$.
\end{pf}

\begin{cor}\label{corPOLsym}
Let $k\geq 3$ and $\ell \geq 2$ be fixed integers.
We can in polynomial time decide if $\kappa_{S,r}(D) \geq \ell$ for any symmetric
digraph, $D$, with $S \subseteq V(D)$, with $|S|=k$ and $ r\in S$.
\end{cor}
\begin{pf}
Let $D$ be any symmetric
digraph with $S \subseteq V(D)$, with $|S|=k$ and $ r\in S$.
Let  ${\cal A} = (A_0,A_1,A_2,\ldots , A_{\ell})$ be a partition of the arcs in $D[S]$. By Theorem~\ref{thmPOLsym} we can decide if
there exist $\ell$ internally disjoint $(S,r)$-trees, $T_1,T_2,\ldots ,T_{\ell}$, with $A(T_i) \cap A[S] = A_i$ for all $i=1,2,\ldots,\ell$.
If such $\ell$ internally disjoint $(S,r)$-trees exist then clearly $\kappa_{S,r}(D) \geq \ell$.
We will do the above for all possible partitions ${\cal A} = (A_0,A_1,A_2,\ldots , A_{\ell})$ and
if we find  $\ell$ internally disjoint $(S,r)$-trees for any such partition then we return ``$\kappa_{S,r}(D) \geq \ell$'' and otherwise we return ``$\kappa_{S,r}(D) < \ell$''.

If  $\kappa_{S,r}(D) \geq \ell$ then we note that we will correctly determine that $\kappa_{S,r}(D) \geq \ell$, when
we consider the correct partition ${\cal A}$, which proves that the above algorithms will always return the correct answer.

Furthermore since the number of partitions, ${\cal A}$, of $A(D[S])$ is bounded by $(\ell+1)^{|A(D[S])|} \leq (\ell+1)^{k^2/2}$
we note that the the above algorithm runs in polynomial time (as $\ell$ and $k$ are considered constants).
\end{pf}

Note that Corollary~\ref{corPOLsym} implies all the polynomial entries in Table~4.

\2

We now turn our attention to the NP-complete cases in Table~4.
Chen, Li, Liu and Mao
\cite{Chen-Li-Liu-Mao} introduced the following problem, which
turned out to be NP-complete.

\2

{\sc CLLM Problem:} Given a tripartite graph $G=(V, E)$ with a
3-partition $(A,B,C)$ such that
$|A|=|B|=|C|=q$, decide whether
there is a partition of $V$ into $q$ disjoint 3-sets $V_1, \dots,
V_q$ such that for every $V_i= \{a_{i_1}, b_{i_2}, c_{i_3}\}$
$a_{i_1} \in A, b_{i_2} \in B, c_{i_3} \in C$ and $G[V_i]$ is connected.

\begin{lem}\label{thm03}\cite{Chen-Li-Liu-Mao} The CLLM Problem is NP-complete.
\end{lem}

In the following theorem we will show the following.
Restricted to symmetric digraphs $D$, for any fixed integer $k\geq
3$, the problem of deciding whether $\kappa_{S, r}(D)\geq \ell~(\ell
\geq 1)$ is NP-complete for $S\subseteq V(D)$ with $|S|=k$ and $r
\in S$.

\begin{thm}\label{thmc}
Let $k\geq 3$ be a fixed integer.
The problem of deciding if a symmetric digraph $D$, with a
$k$-subset $S$ of $V(D)$ with $r\in S$ satisfies $\kappa_{S, r}(D)\geq \ell$ ($\ell$ is part of the input), is
NP-complete.
\end{thm}
\begin{pf}
It is easy to see that this problem is in NP. Let $G$ be a
tripartite graph with 3-partition $(A,B,C)$ such that
$|A|=|B|=|C|=q$. We will construct
a symmetric digraph $D$, a $k$-subset $S\subseteq V(D)$ with $r\in
S$ and an integer $\ell$ such that there are $\ell$ internally
disjoint $(S, r)$-trees in $D$ if and only if $G$ is a positive
instance of the CLLM Problem.

Let $D$ be obtained from $G$ by replacing every edge with a $2$-cycle and
adding the vertices $S=\{r,s_1,s_2,\ldots,s_{k-1}\}$ and all arcs between
$r$ and $A$ and all arcs between $s_1$ and $B$ and all arcs between $\{s_2,s_3,\ldots,s_{k-1}\}$ and $C$.
This completes the construction of $D$.
Note that the construction of $D$ is equivalent to a construction in \cite{Sun-Gutin-Yeo-Zhang} and clearly
can be done in polynomial time.

First consider the case when there are  $\ell$ internally disjoint $(S, r)$-trees in $D$, say
$T_i~(1\leq i\leq \ell)$. Each tree must clearly contain at least one vertex from $A$ (connected to $r$ in the tree),
at least one vertex from $B$ (connected to $s_1$ in the tree) and
at least one vertex from $C$ (connected to $s_2$ in the tree). As $|A|=|B|=|C|=q$ and the
trees are internally disjoint, we note that every tree, $T_i$, contains exactly one vertex from $A$, say $a_i$, and one vertex from $B$,
say $b_i$,  and exactly one vertex from $C$, say $c_i$. However now $G[{a_i,b_i,c_i}]$ is connected for all $i=1,2,\ldots,q$,
and $G$ is a positive instance of CLLM.

Conversely if $G$ is a positive instance of CLLM, then there is a partition of $V(G)$ into $q=\ell$ disjoint sets $V_1,
V_2, \cdots, V_q$ each having three vertices, such that for every
$V_i= \{a_{i_1}, b_{i_2}, c_{i_3}\}$ we have $a_{i_1} \in A$,
 $b_{i_2} \in B$ and $c_{i_3} \in C$, and $G[V_i]$ is connected.
Let $T_i$ be a $(S,r)$-tree with vertex-set $V_i \cup S$.
Note that all $T_i$ are internally disjoint $(S,r)$-trees in $D$.
By the above argument and Lemma~\ref{thm03}, we are done.
\end{pf}

Theorem~\ref{thmc} together with the below theorem implies all the NP-completness results in Table~4.

\begin{thm}\label{thmLfree}
Let $\ell \geq 2$ be a fixed integer.
The problem of deciding if a symmetric digraph $D$, with a
$S \subseteq V(D)$ and $r\in S$ satisfies $\kappa_{S, r}(D)\geq \ell$ ($k=|S|$ is part of the input), is NP-complete.
\end{thm}
\begin{pf}
  We will reduce from the NP-hard problem of $2$-coloring a hypergraph (see \cite{Lovasz}).
That is, we are given a hypergraph, $H$, with vertex-set $V(H)$ and edge-set $E(H)$, and want to
determine if we can $2$-colour the vertices $V(H)$ such that every hyperedge in $E(H)$ contains
vertices of both colours. This problem is known to be NP-hard (see \cite{Lovasz}).

  Define a symmetric digraph, $D$, as follows.
Let $U = \{u_1,u_2,\ldots, u_{\ell-2} \}$ and
let $V(D)= V(H) \cup E(H) \cup U \cup \{r\}$ and let the arc-set of $D$ be defined as follows.

\[
\begin{array}{rcl}
A(D) & = & \{xe ,ex \; | \; x\in V(H), \mbox{ } e \in e(H) \mbox{ and } x \in V(e)\} \\
     &   &  \cup \; \; \{r u_i, u_i r, u_i e, e u_i \; | \; u_i \in U \mbox{ and } e \in E(H) \} \\
     &   &  \cup \; \;  \{r x, x r \; | \; x \in V(H)\} \\
\end{array}
\]

Let $S = E(H) \cup \{r\}$. This completes the construction of $D$, $S$ and $r$.
We will show  that $\kappa_{S, r}(D)\geq \ell$ if and only if $H$ is 2-colourable, which
will complete the proof.

First assume that $H$ is $2$-colourable and let $R$ be the red vertices in $H$ and $B$ be the blue vertices in $H$ in
a proper $2$-colouring of $H$.
Let $T_i$ contain the arc $r u_i$ and all arcs from $u_i$ to $S \setminus \{r\}$ for $i=1,2,\ldots, \ell-2$.
Let $T_{\ell-1}$ contain all arcs from $r$ to $R$ and for each edge $e \in E(H)$ we add an arc from $R$ to $e$ in $D$ to $T_{\ell-1}$
(this is possible as every edge in $H$ contains a red vertex).
Analogously, let $T_{\ell}$ contain all arcs from $r$ to $B$ and for each edge $e \in E(H)$ we add an arc from $B$ to $e$ in $D$ to $T_{\ell}$
(again, this is possible as every edge in $H$ contains a blue vertex).  We now note that $T_1,T_2,\ldots , T_{\ell}$ are internally disjoint $(S,r)$-trees in $D$,
so $\kappa_{S, r}(D)\geq \ell$.

Conversely assume that $\kappa_{S, r}(D)\geq \ell$ and let $T_1',T_2',\ldots,T_{\ell}'$ be $\ell$ internally disjoint $(S,r)$-trees in $D$.
At least two of these trees contain no vertex from $U$ (as $|U| = \ell-2$). Without loss of generality assume that
$T_1'$ and $T_2'$ don't contain any vertex from $U$. Let $B'$ be all out-neighbours of $r$ in $T_1'$ (i.e. $B' = N_{T_1'}^+(r)$)
and let $R'$ be all out-neighbours of $r$ in $T_2'$  (i.e.  $R' = N_{T_2'}^+(r)$). For every $e \in E(H)$ it has an arc into it in
$T_1'$ and an arc into it in $T_2'$ which implies that in $H$ the edge $e$ contains a vertex from $B'$ and a vertex from $R'$.
Therefore $H$ is $2$-colourable (any vertex in $H$ that is not in either $R'$ or $B'$ can be assigned randomly to either $B'$ or $R'$).
This completes the proof.
\end{pf}

\section{Proof of Theorem~\ref{LambdaEulerian} (Eulerian digraphs)}\label{sec:proofL}

We will in this section give a proof for Theorem~\ref{LambdaEulerian}. However we first need the following Theorem by Bang-Jensen, Frank and Jackson.

\begin{thm} \cite{JBJ-Frank-Jackson} \label{ThmEulerianJBJ}
Let $k \geq 1$ and let $D=(V,A)$ be a directed multigraph with a special vertex $z$.
Let $T' = \{ x \; | \; x \in V \setminus \{z\} \mbox{ and } d^-(x) < d^+(x) \}$.
If $\lambda(z,x) \geq k$ for every $x \in T'$, then there exists a family ${\cal F}$ of $k$ arc-disjoint out-trees rooted at $z$ so that every vertex
$x \in V$ belongs to at least $\min\{k,\lambda(z,x)\}$ members of ${\cal F}$.
\end{thm}

In the case when $D$ is Eulerian, then $d^+(x)=d^-(x)$ for all $x \in V(D)$ and therefore $T'=\emptyset$ in the above theorem.  Therefore the following
corollary holds.

\begin{cor} \label{CorEulerianJBJ}
Let $k \geq 1$ and let $D=(V,A)$ be an Eulerian digraph with a special vertex $z$.
Then there exists a family ${\cal F}$ of $k$ arc-disjoint out-trees rooted at $z$ so that every vertex
$x \in V$ belongs to at least $\min\{k,\lambda(z,x)\}$ members of ${\cal F}$.
\end{cor}

We are now in a position to prove Theorem~\ref{LambdaEulerian}. Recall the theorem.

\2

\noindent
{\bf Theorem~\ref{LambdaEulerian}:} {\em   If $D$ is an Eulerian digraph and $S \subseteq V(D)$ and $r \in S$, then
$\lambda_{S,r}(D) \geq \ell$ if and only if $\lambda_D(r,s) \geq \ell$ for all $s \in S \setminus \{r\}$.}

\2

\begin{pf}
Let $D$ be an Eulerian digraph and let $S \subseteq V(D)$ and $r \in S$ be arbitrary.
First assume that $\lambda_{S,r}(D) \geq \ell$. This implies that $\lambda_D(r,s) \geq \ell$ for all $s \in  S \setminus \{r\}$
as there is a path from $r$ to $s$ in each of the $\ell$ arc-disjoint $(S,r)$-trees and these $\ell$ paths are therefore also arc-disjoint.

Now assume that $\lambda_D(r,s) \geq \ell$ for all $s \in S \setminus \{r\}$. By Corollary~\ref{CorEulerianJBJ} there exists
a family ${\cal F}$ of $\ell$ arc-disjoint out-trees rooted at $r$ so that every vertex
$x \in V$ belongs to at least $\min\{\ell,\lambda(r,x)\}$ members of ${\cal F}$.
As $\lambda(r,s) \geq \ell$ for all $s \in S \setminus \{r\}$ we note that every vertex in $S$ belongs to all $\ell$
out-trees in ${\cal F}$.  Therefore $\lambda_{S,r}(D) \geq \ell$, which completes the proof of the theorem.
\end{pf}

As one can determine $\lambda_D(r,s)$ in polynomial time for any $r$ and $s$ in $D$ we note that
Theorem~\ref{LambdaEulerian} implies all the entries in Tables~3 and 5.

In the proof of Theorem~\ref{thmPOLsym} we saw how the fact that the $k$-linkage problem in symmetric graphs is polynomial
was used to prove the polynomial entries in Table~4. It is maybe therefore slightly surprising that even though we prove that
all entries in Table~5 are polynomial, this could not be proved using a similar linkage result, due to the following.

\begin{thm} \cite{Ibaraki-Poljak}
The weak $k$-linkage problem is NP-hard in Eulerian digraphs (where $k$ is part of the input).
\end{thm}

\section{Entries of Table~6 (Eulerian digraphs)}

We will in this section prove the NP-completeness results given in Table~6.
In our argument, we will use the constructions from the proof for the following theorem by Fortune, Hopcroft and Wyllie. And in fact, we will use Figures 10.1-10.3 in the arugment of Theorem~10.2.1 of \cite{Bang-Jensen-Gutin}.

\begin{thm}\label{2link}\cite{Fortune-Hopcroft-Wyl}
The 2-linkage problem is NP-complete.
\end{thm}

\begin{figure}[!hbpt]
\begin{center}
\includegraphics[scale=1.0]{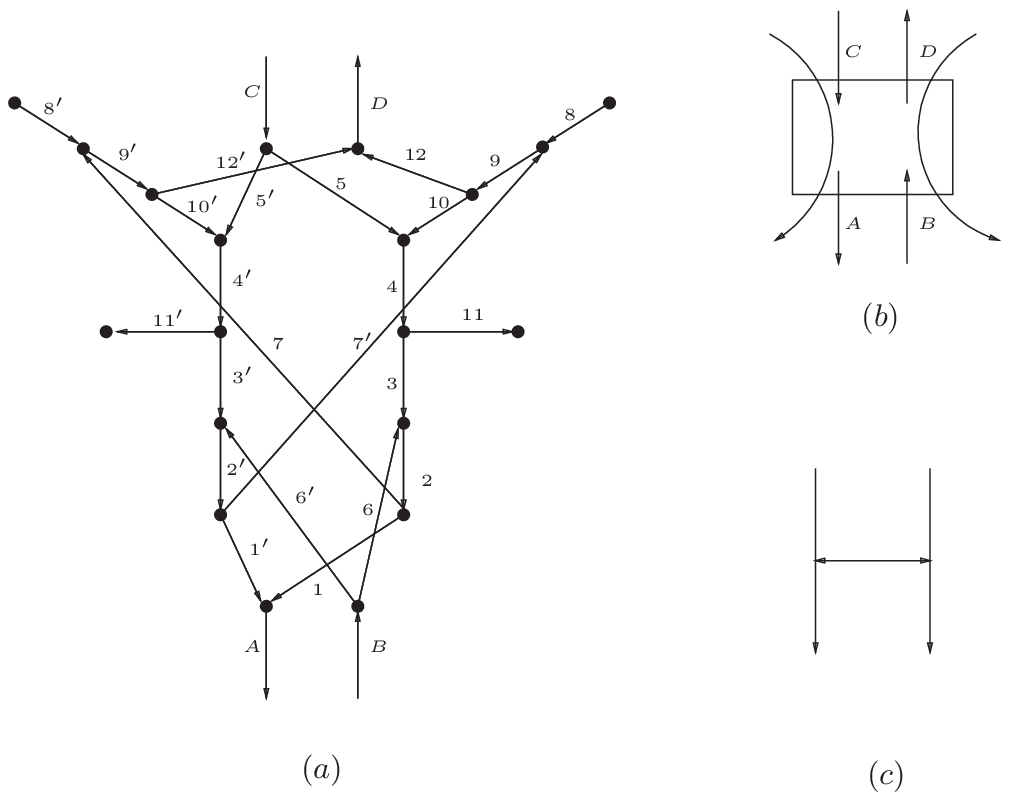}
\end{center}
\caption{A switch $(a)$ and its schematic pictures $(b)$ and $(c)$\cite{Bang-Jensen-Gutin, Fortune-Hopcroft-Wyl}.}\label{figure01}
\end{figure}

In their argument, Fortune, Hopcroft and Wyllie used the concept of {\em switch} which is shown in Figure~\ref{figure01}$(a)$. Note that in $(c)$ the two vertical arcs correspond to the paths
(8, 9, 10, 4, 11), respectively, (8', 9', 10', 4', 11'). For convenience, we label the
arcs, rather than the vertices in this figure.

By the definition of a switch, the following lemma holds:

\begin{lem}\label{2link-1}\cite{Fortune-Hopcroft-Wyl}
Consider the digraph $S$ shown in Figure~\ref{figure01}$(a)$. Suppose there are two vertex-disjoint paths $P, Q$ passing through $S$ such that $P$ leaves $S$ at $A$ and $Q$ enters $S$ at $B$. Then $P$ must enter $S$ at $C$ and $Q$ must leave $S$ at $D$. Furthermore, there exists exactly one more path $R$ passing through
$S$ which is disjoint from $P, Q$ and this is either (8, 9, 10, 4, 11) or (8', 9', 10', 4', 11'),
depending on the actual routing of $P$.
\end{lem}

\begin{figure}[!hbpt]
\begin{center}
\includegraphics[scale=1.0]{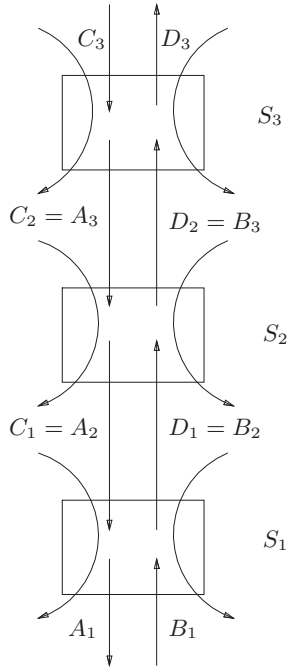}
\end{center}
\caption{Stacking three switches on top of each other \cite{Bang-Jensen-Gutin, Fortune-Hopcroft-Wyl}.}\label{figure02}
\end{figure}

As shown in \cite{Bang-Jensen-Gutin, Fortune-Hopcroft-Wyl}, we can stack arbitrarily many switches on top of each other and still have the conclusion on Lemma~\ref{2link-1} holding for each switch. The way we stack is simply by identifying the $C$ and $D$ arcs of one switch with the $A$ and $B$ arcs of the next (see Figure~\ref{figure02}). 

Now we can prove the NP-completeness of 2-linkage problem for Eulerian digraphs.

\begin{thm}\label{NP-eulerian}
The 2-linkage problem restricted to Eulerian digraphs is NP-complete.
\end{thm}
\begin{pf}
The reduction of the argument for Theorem~\ref{2link} is from 3-SAT problem. Let $\mathcal{F}=C_1 \wedge C_2 \wedge \cdots \wedge C_r$ be an instance of 3-SAT with varibles $x_1, x_2, \dots, x_k$ and clauses $C_1,C_2,\ldots,C_r$.
For each variable $x_i$ we let $H_i$ be a digraph consisting of two internally disjoint $(z_i, w_i)$-paths of length $r$ (the number of clauses in $F$). See Figure~\ref{figure03} for an illustration. We associate one of these paths with the literal $x_i$ and the other with $\overline{x}_i$.

\begin{figure}[!hbpt]
\begin{center}
\includegraphics[scale=1.0]{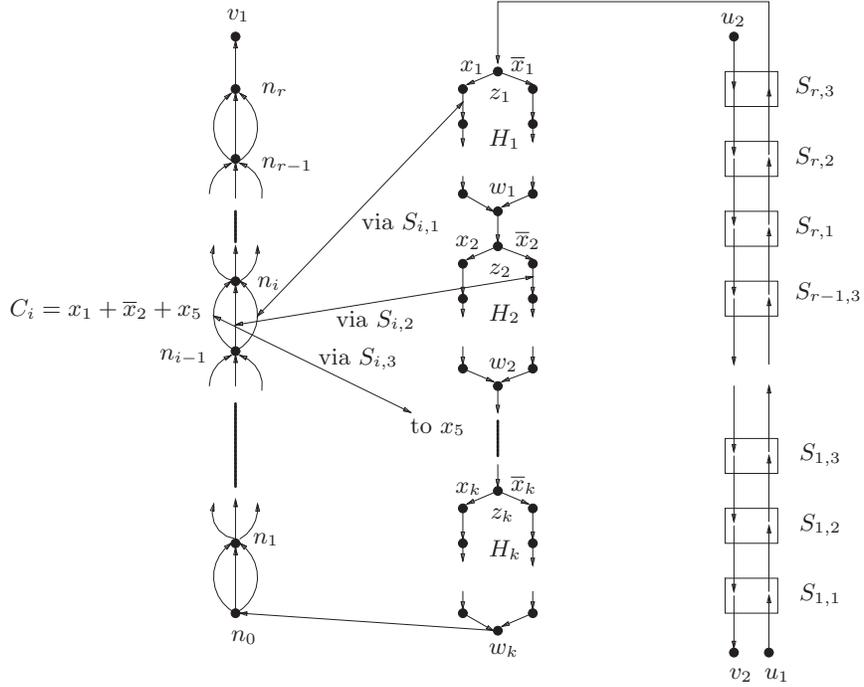}
\end{center}
\caption{ A schematic picture of the digraph $D[\mathcal{F}]$\cite{Bang-Jensen-Gutin, Fortune-Hopcroft-Wyl}.}\label{figure03}
\end{figure}

We will now construct the digraph $D[\mathcal{F}]$ as follows (see Figure~\ref{figure03} and e.g. \cite{Bang-Jensen-Gutin}).
Firstly, we form a chain $H_1\rightarrow H_2\rightarrow \dots \rightarrow H_k$ on the subdigraphs corresponding to each variable (see the middle of the figure, $H_i$ corresponds
to the variable $x_i$). Secondly, with each clause $C_i$ we associate three switches, one for
each literal it contains, and they are denoted by $S_{i, 1}, S_{i, 2}, S_{i, 3}$, respectively;
we then stack these switches in the order $S_{1, 1}, S_{1, 2}, S_{1, 3}, S_{2,1}, \dots, S_{r, 1}, S_{r,2}, S_{r,3}$ as shown in the right part of the figure.

Thirdly, we create a path $n_0, n_1,\ldots,n_r$, such that there are three arcs from $n_{i-1}$ to $n_i$ for all $i=1,2,\ldots,r$, each corresponding to a literal in $C_i$.
And if $x_a$ is the $b$'th literal in $C_i$, then we substitute the $b$'th arc from $n_{i-1}$ to $n_i$ with the "left" path of $S_{i,b}$ (that is,
the path with arcs $8',9',10',4',11'$ in $S_{i,b}$) and we substitute a (private) arc of $H_a$, such that the arc is taken from the path which corresponds to $x_a$ if the literal is $x_a$ and from the path which corresponds to $\overline{x}_a$ if the literal is $\overline{x}_a$,
with the "right" path of  $S_{i,b}$ (that is,
the path with arcs $8,9,10,4,11$ in $S_{i,b}$).
For example in Figure~\ref{figure03} this would imply that the right-most arc from $n_{i-1}$ to $n_i$ actually is the path
$8',9',10',4',11'$ in $S_{i,1}$ (as $x_1$ is the first literal in $C_i$) and the arc indicated in $H_1$ is actually the path
$8,9,10,4,11$ in $S_{i,1}$.
In this way one can show that only one of the paths $8',9',10',4',11'$  and $8,9,10,4,11$ in each $S_{i,b}$ will ever be used in
a solution to the $2$-linkage problem. In other words, The double arcs which point to two different arcs indicate that only one of these
arcs (which are actually paths) can be used in a solution.

Finally, we join the $D$ arc of the switch $S_{r,3}$ to the vertex $z_1$ of $H_1$, add an arc from $w_k$ in $H_k$ to $n_0$ and choose vertices $u_1, u_2, v_1, v_2$ as shown in the figure.

We now construct an Eulerian digraph $D'$ from $D[\mathcal{F}]$ as follows:
Firstly, for each switch $S_{i,j}$, we duplicate the arcs $A, B, C, D, 9', 9, 4', 4, 2, 2'$ (note that in this procedure there are two parallel arcs between $S_{r,3}$ and the vertex $z_1$).
Secondly, We duplicate the arc $w_kn_0$, and the arc $w_iz_{i+1}$ for each $1\leq i\leq k-1$, respectively; we also duplicate the arc $n_rv_1$ twice, and add the arc $v_1n_0$. Observe that now each $u_i$ has exactly two arcs into $U=V(D[\mathcal{F}])\setminus \{u_1, u_2, v_1, v_2\}$ and has no arc from it, $v_1$ has three arcs from $U$ and one arc into it, $v_2$ has two arcs from $U$ and has no arc into it. Finally, we add a new vertex $x$ and the following new arcs: add the arcs $v_ix$ and $xu_j$ twice for each $1\leq i, j\leq 2$. To avoid parallel arcs, we could subdivide each new arc.

It is not difficult to show that $D[\mathcal{F}]$ contains a pair of vertex-disjoint $(u_1, v_1)$-, $(u_2, v_2)$-paths if and only if $D'$ contains a pair of vertex-disjoint $(u_1, v_1)$-, $(u_2, v_2)$-paths. Indeed, if $D[\mathcal{F}]$ contains a pair of vertex-disjoint $(u_1, v_1)$-, $(u_2, v_2)$-paths, then clearly these two paths are also vertex-disjoint $(u_1, v_1)$-, $(u_2, v_2)$-paths in $D'$. For the other direction, let $P$ and $Q$ be vertex-disjoint $(u_1, v_1)$-, $(u_2, v_2)$-paths in $D'$, then clearly they do not contain arcs incident to $x$ and the arc $v_1n_0$. If they use other new arcs, then we can replace them by the corresponding parallel original arcs in $D[\mathcal{F}]$, and obtain desired paths in $D[\mathcal{F}]$. It was proved in \cite{Fortune-Hopcroft-Wyl} that $D[\mathcal{F}]$ contains a pair of vertex-disjoint $(u_1, v_1)$-, $(u_2, v_2)$-paths if and only if $\mathcal{F}$ is satisfiable. Hence, the Eulerian digraph $D'$ contains a pair of vertex-disjoint $(u_1, v_1)$-, $(u_2, v_2)$-paths if and only if $\mathcal{F}$ is satisfiable, and therefore the 2-linkage problem for Eulerian digraphs is NP-complete.
\end{pf}

Using Theorem~\ref{NP-eulerian}, we can prove the following result, which completes all the entries in Table~6.

\begin{thm}\label{tree-eulerian}
Let $k\geq 3$ and $\ell\geq 2$ be fixed integers. Let $D$ be an Eulerian
digraph and $S \subseteq V(D)$ with $|S|=k$ and $r \in S$. Then deciding whether $\kappa_{S, r}(D)\ge \ell$ is NP-complete.
\end{thm}
\begin{pf}
It is not difficult to see that the problem belongs to NP.
We will show that the problem is NP-hard by reducing from $2$-linkage in Eulerian digraphs.
Let $H$ be an Eulerian digraph and let $s_1,s_2,t_1,t_2$ be distinct vertices in $H$.
We now produce a new Eulerian digraph $D^*$ with $V(D^*) = V(H) \cup V \cup U \cup \{r\}$,
where $V=\{v_1,v_2, \ldots v_{\ell-2}\}$ and $U=\{u_1, u_2, \ldots, u_{k-1}\}$ (here let $S=\{r, u_1, u_2, \ldots, u_{k-1}\}$).
Furthermore let the arcs set of $D^*$ be defined as follows.

\begin{figure}
\begin{center}
\tikzstyle{vertexB}=[circle,draw, minimum size=15pt, scale=0.7, inner sep=0.1pt]
\tikzstyle{vertexL}=[circle,draw, minimum size=15pt, scale=0.9, inner sep=0.1pt]
\begin{tikzpicture}[scale=0.20]
\node (r) at (-2,6) [vertexL] {$r$};
\node (s1) at (8,4) [vertexL] {$s_1$};
\node (s2) at (8,8) [vertexL] {$s_2$};
\node (t1) at (28,2) [vertexL] {$t_1$};
\node (t2) at (28,10) [vertexL] {$t_2$};
\draw (5,0) rectangle (31,12);

\node (u1) at (40,4) [vertexL] {$u_1$};
\node (u2) at (40,8) [vertexL] {$u_2$};
\draw (38,0) rectangle (42,12);

\node (u3) at (40,17) [vertexL] {$u_3$};
\node (u4) at (40,20) [vertexL] {$u_4$};
\node (uk) at (40,26) [vertexB] {$u_{k-1}$};
\draw (38,15) rectangle (42,28);


\node (v1) at (15,24) [vertexL] {$v_1$};
\node (v2) at (18,24) [vertexL] {$v_2$};
\node (vl) at (24,24) [vertexB] {$v_{\ell-2}$};
\draw (12,21) rectangle (27,27);
\node () at (19.5,28) [scale=1.1] {$V$};

\node () at (20,3) [scale=1.5] {$H$};
\node () at (21,26) [scale=0.8] {$\cdots$};
\node () at (40,23) [scale=0.8] {$\vdots$};


\draw [->, line width=0.02cm] (s1) to (u2);
\draw [->, line width=0.02cm] (s2) to (u1);
\draw [->, line width=0.02cm] (r) to (s2);
\draw [->, line width=0.02cm] (r) to (s1);

\draw [<->, line width=0.04cm] (r) to [out=90, in=200] (12,23);
\draw [<->, line width=0.04cm] (t1) to (u1);
\draw [<->, line width=0.04cm] (t2) to (u2);
\draw [<->, line width=0.04cm] (40,12) to (40,15);
\draw [<->, line width=0.04cm] (27,24) to (38,24);
\draw [<->, line width=0.04cm] (27,21) to (38,12);

\draw [line width=0.02cm] (u2) to [out=130, in=0] (26,16);
\draw [line width=0.02cm] (26,16) to (10,16);
\draw [->, line width=0.02cm] (10,16) to [out=180, in=70] (r);

\draw [line width=0.02cm] (u1) to [out=230, in=0] (26,-4);
\draw [line width=0.02cm] (26,-4) to (10,-4);
\draw [->, line width=0.02cm] (10,-4) to [out=180, in=290] (r);

\end{tikzpicture}
\end{center}
\caption{Illustration of $D^*$ in the proof of Theorem~\ref{tree-eulerian}.}
\end{figure}
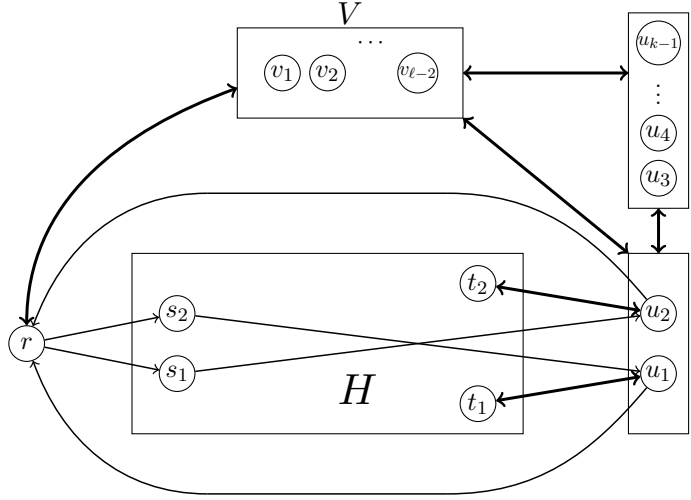

\[
\begin{array}{rcl}
 A(D^*) & = & A(H) \cup \{rs_1, rs_2, t_1 u_1, u_1t_1, t_2u_2, u_2t_2, s_1u_2, s_2u_1, u_1r, u_2r\} \\
        &   & \cup \; \{ rv, vr, vu, uv \; | \; v \in V \mbox{ and } u \in U \} \\
        &   & \cup \; \{u_i u_j, u_j u_i \; | \; i=1,2 \mbox{ and } j=3,4,\ldots, k-1\} \\
\end{array}
\]

Note that $D^*$ is Eulerian and let $S=\{r\} \cup U$. We will show that $\kappa_{S, r}(D^*)\ge \ell$ if and only if
there exist two vertex-disjoint paths, $P_1$ and $P_2$, such that $P_i$ is an $(s_i,t_i)$-path, for $i=1,2$.
This will complete the proof.

First assume that there exist two vertex-disjoint paths, $P_1$ and $P_2$, such that $P_i$ is an $(s_i,t_i)$-path, for $i=1,2$.
Add the arcs, $r s_1, s_1 u_2, t_1 u_1$ and all arcs from $u_1$ to $\{u_3,u_4,\ldots,u_{k-1}\}$ to $P_1$ and call the resulting
$(S,r)$-tree for $T_{\ell -1}$.  Analogously, add the arcs, $r s_2, s_2 u_1, t_2 u_2$ and all arcs from $u_2$ to $\{u_3,u_4,\ldots,u_{k-1}\}$ to $P_2$
and call the resulting $(S,r)$-tree for $T_{\ell}$. Finally let $T_i$ be the $(S,r)$-tree containing the arc $r v_i$ and all arcs from $v_i$ to $U$, for
each $i=1,2,3,\ldots,\ell-2$.  The $(S,r)$-trees, $T_1,T_2,\ldots,T_{\ell}$ are now internally disjoint, which implies that  $\kappa_{S, r}(D^*)\ge \ell$
as desired.

Conversely assume that $\kappa_{S, r}(D^*)\ge \ell$ and let $T_1$ and $T_2$ be two $(S,r)$-trees in $D^*$ that do not use any of the vertices in $V$ (which exist as $|V|=\ell-2$).  Without loss of generality assume that $r s_1 \in A(T_1)$ and $rs_2 \in A(T_2)$. As the arcs into $u_1$ in $T_1$ and $T_2$
is either $s_2u_1$  or $t_1 u_1$ we note that $s_2 u_1 \in A(T_2)$ (as $s_2 \in V(T_2)$) and $t_1 u_1 \in A(T_1)$.
Analogously,  $s_1 u_2 \in A(T_1)$ and $t_2 u_2 \in A(T_2)$. This implies that there is an $(s_1,t_1)$-path in $T_1$ and an
$(s_2,t_2)$-path in $T_2$, which are vertex-disjoint.  This competes the proof of the fact that
deciding whether $\kappa_{S, r}(D)\ge \ell$ is NP-complete.
\end{pf}

\section{Sharp bounds and characterizations for $\kappa_k(D)$ and $\lambda_k(D)$}\label{sec:bounds}

The following proposition can be verified using definitions of
$\kappa_k(D)$ and $\lambda_{k}(D)$.

\begin{pro}\label{proa}
Let $D$ be a digraph of order $n$, and let $k\ge 2$ be an integer.
Then
\begin{equation}\label{pro1}
\lambda_{k+1}(D)\leq \lambda_{k}(D) \mbox{ for every } k\le n-1
\end{equation}
\begin{equation}\label{pro2}
\kappa_k(D')\leq \kappa_k(D), \lambda_k(D')\leq \lambda_k(D) \mbox{
where $D'$ is a spanning subdigraph of $D$}
\end{equation}
\begin{equation}\label{pro3}
\kappa_k(D)\leq \lambda_k(D) \leq \min\{\delta^+(D), \delta^-(D)\}
\end{equation}
\end{pro}

Note that Proposition \ref{proa}(\ref{pro1}) may not hold for $\kappa_{k}$, that is, $\kappa_n(D)\leq \kappa_{n-1}(D)\leq \ldots \kappa_3(D)\leq \kappa_2(D)=\kappa(D)$ may not be true.
Consider the following example: Let $D_1$ and $D_2$ be two copies of the complete digraphs $\overleftrightarrow{K}_{t}~(t\geq 4)$,
and let $D$ be a digraph obtained from $D_1$ and $D_2$ by identifying one vertex in each of them.
Clearly, $D$ is a strong connected digraph with a cut vertex, say $u$. For $2\leq k\leq 2t-2$,
let $S$ be a subset of $V(D)\setminus \{u\}$ with $|S|=k$.
Since each $(S,r)$-tree must contain $u$, we have $\kappa_k(D)\leq 1$, furthermore, we deduce that $\kappa_k(D)= 1$ for $2\leq k\leq 2t-2$.
Let $G_i$ be the underlying undirected graph of $D_i$ for $1\leq i\leq 2$. Each $G_i$ contains $\lfloor \frac{t}{2}\rfloor$ edge-disjoint spanning trees, say $T_{i,j}~(1\leq j\leq \lfloor \frac{t}{2}\rfloor)$, since $G_i$ is a complete graph of order $t$, where $1\leq i\leq 2$. Now in $D$, let $D_j$ be a subdigraph obtained from the tree $T_j$ which is the union of $T_{1,j}$ and $T_{2,j}$ by replacing each edge with the corresponding arcs of both directions. Clearly, these subdigraphs are strong spanning subdigraphs and pair-wise arc-disjoint. From these subdigraphs we can find at least $\lfloor \frac{t}{2}\rfloor$ arc-disjoint out-branchings rooted at each vertex of $D$ by the fact that every strong digraph has an out- and in-branching rooted at any vertex. Hence, $\kappa_{2t-1}(D)\geq \lfloor \frac{t}{2}\rfloor>1=\kappa_k(D)$ for $2\leq k\leq 2t-2$.


The following result concerning the exact values of $\kappa_k(\overleftrightarrow{K}_n)$ and
$\lambda_k(\overleftrightarrow{K}_n)$ will be used in the proof of one of our main results in this section.

\begin{lem}\label{thm1} For $2\leq k\leq n$, we have
$$\kappa_k(\overleftrightarrow{K}_n)=\lambda_k(\overleftrightarrow{K}_n)=n-1.$$
\end{lem}
\begin{pf} Let $D=\overleftrightarrow{K}_n$.
By Proposition~\ref{proa}(\ref{pro3}), we
have $\kappa_k(D)\leq \delta^+(D)= n-1$ and $\lambda_k(D)\leq \delta^+(D)= n-1$, so we just need to show that there exists
$n-1$ internally disjoint $(S,r)$-trees in $D$, for any $S \subseteq V(D)$ with $|S|=k$ and $r \in S$.
Let $V(D)=\{u_1,u_2,\ldots,u_n\}$ and without loss of generality that $S=\{u_1,u_2,\ldots,u_k\}$ and $r=u_1$.

For $i=2,3,\ldots, k$ let $T_i$ be the $(S,r)$-tree containing the arc $r u_i$ and all
arcs from $u_i$ to $S \setminus \{r,u_i\}$.
Note that thess $k-1$ $(S,r)$-trees are arc-disjoint and internally disjoint (as they all lie completely within $S$).

For $i=k+1,k+2,\ldots,n$ let $T_i$ be the $(S,r)$-tree containing the arc $r u_i$ and all
arcs from $u_i$ to $S \setminus \{r\}$. Now all the $(S,r)$-trees, $T_2,T_3,\ldots,T_n$
are arc-disjoint and internally disjoint.
Therefore $\kappa_k(D) = \lambda_k(D) = n-1$.
\end{pf}

\begin{thm}\label{thma}
Let $D$ be a strong digraph of order $n$, and let $k\ge 2$ be an
integer. Then
\begin{equation}\label{pro4}
1\leq \kappa_k(D)\leq n-1
\end{equation}
\begin{equation}\label{pro5}
1\leq \lambda_k(D)\leq n-1
\end{equation} Moreover, all bounds are sharp, and
the upper bounds hold if and only if $D\cong
\overleftrightarrow{K}_n$.
\end{thm}
\begin{pf}
The upper bounds hold by Proposition~\ref{proa}(\ref{pro2}) and
(\ref{pro3}), and lower bounds hold by the fact that every strong digraph has an out- and in-branching rooted at any vertex. For the
sharpness of the lower bound, a cycle is our desired digraph. If $D$
is not equal to $\overleftrightarrow{K}_n$ then $\delta^+(D) \leq
n-2$ and by Proposition~\ref{proa} we note that $\kappa_k(D) \leq
\delta^+(D) \leq n-2$ and $\lambda_k(D) \leq \delta^+(D) \leq n-2$.
Therefore, by Lemma~\ref{thm1}, the upper bounds hold if and only if
$D\cong \overleftrightarrow{K}_n$.
\end{pf}

The following sharp bounds concerning the relation between
$\kappa_k(D)$~(respectively, $\lambda_k(D))$ and $\kappa(D)$~(respectively, $\lambda(D))$ improve
those of Proposition~\ref{proa}(\ref{pro3}).

\begin{thm}\label{thmb}
Let $2\leq k\leq n$ be an integer. The following assertions hold:
\begin{description}
\item[(i)]~$\kappa_k(D)\leq \kappa(D)$ when $n\ge \kappa(D)+k$.
\item[(ii)]~$\lambda_k(D)\leq \lambda(D).$
\end{description}
Moreover, both bounds are sharp.
\end{thm}
\begin{pf} \noindent{\bf Part $(i)$}.
For $k=2$, we have $\kappa_2(D)= \kappa(D)$ by definition. In the
following argument we consider the case of $k\ge 3$.

If $\kappa(D)=0$, then $D$ is not strong and $\kappa_k(D)=0$, as can be seen by letting $r,x \in S$ be chosen such that there is
no $(r,x)$-path in $D$.
If $\kappa(D)=n-1$, then we have $\kappa_k(D)\leq n-1=\kappa(D)$ by
(\ref{pro4}) in Theorem~\ref{thma}, so we may assume that $ 1 \leq \kappa(D) \leq n-2$. There now exists a $\kappa(D)$-vertex
cut, say $Q$, for two vertices $u,v$ in $D$ such that there is no
$(u,v)$-path in $D-Q$. Let $S=\{u,v\}\cup S'$ where $S'\subseteq
V(D)\setminus (Q\cup \{u,v\})$ and $|S'|=k-2$. Observe that in each
$(S, u)$-tree, the $u-v$ path must contain a vertex in $Q$. By the
definition of $\kappa_{S,r}(D)$ and $\kappa_k(D)$, we have
$\kappa_k(D)\leq \kappa_{S,r}(D)\leq |Q|=\kappa(D)$.

For the sharpness of the bound in $(i)$, 
consider the following digraph $D$. Let $D$ be a symmetric digraph
whose underlying undirected graph is $K_{k}\bigvee
\overline{K}_{n-k}$~($n\geq 3k$), i.e. the graph obtained from
disjoint graphs $K_{k}$ and $\overline{K}_{n-k}$ by adding all edges
between the vertices in $K_{k}$ and $\overline{K}_{n-k}$. 

Let $V(D)=W\cup U$, where $W=V(K_k)=\{w_i\mid 1\leq i\leq k\}$ and
$U=V(\overline{K}_{n-k})=\{u_j\mid 1\leq j\leq n-k\}$. Note that $n-k\geq 2k$ since $n\geq 3k$. Let $S$ be
any $k$-subset of vertices of $V(D)$ such that $|S\cap U|=s$ ($s\leq
k$) and $|S\cap W|=k-s$. Without loss of generality, let $w_i\in S$
for $1\leq i\leq k-s$ and $u_j\in S$ for $1\leq j\leq s$. For $1\leq
i\leq k-s$, let 
$T_i$ be a tree with edge set
$$\{w_iu_1, w_iu_2, \dots , w_iu_s, u_{k+i}w_1, u_{k+i}w_2, \dots ,
u_{k+i}w_{k-s}\}.$$ For $k-s+1\leq j\leq k$, let 
$T_j$ be a tree with edge set $$\{w_ju_1, w_ju_2, \dots , w_ju_s,
w_jw_1, w_jw_2, \dots , w_jw_{k-s}\}.$$ It is not hard to obtain an
$(S, r)$-tree $D_i$ from $T_i$ by adding appropriate directions to
edges of $T_i$ for any $r\in S$. Observe that $\{D_i\mid 1\leq i\leq
k-s\}\cup \{D_j\mid k-s+1\leq j\leq k\}$ is a set of $k$ internally
disjoint $(S, r)$-trees, so $\kappa_{S, r}(D)\geq k$, and then
$\kappa_k(D)\geq k$. Combining this with the bound that
$\kappa_k(D)\leq \kappa(D)$ and the fact that $\kappa(D)\leq
\min\{\delta^+(D), \delta^-(D)\}=k$, we have $\kappa_k(D)=
\kappa(D)=k$.

\noindent{\bf Part $(ii)$}. The bound in $(ii)$ is from
Proposition~\ref{proa}(\ref{pro1}) and the fact that $\lambda_2(D)=\lambda(D)$. For the sharpness of this bound,
just consider the above example. Clearly, $\{D_i\mid 1\leq i\leq
k-s\}\cup \{D_j\mid k-s+1\leq j\leq k\}$ is also a set of $k$
arc-disjoint $(S, r)$-trees for any $r\in S$, so $\lambda_{S,
r}(D)\geq k$, and then $\lambda_k(D)\geq k$. Combining this with the
bound that $\lambda_k(D)\leq \lambda(D)=k$, we can get
$\lambda_k(D)= \lambda(D)=k$. This completes the proof.
\end{pf}

Note that the condition ``$n\ge \kappa(D)+k$" in Theorem \ref{thmb}$(i)$ cannot be removed. Consider the example before Lemma~\ref{thm1}, we have $n=2t-1< 2t=\kappa(D)+k$ when $k=n$, but now $\kappa_n(D)>\kappa(D)$.

Given a graph parameter $f(G)$, the Nordhaus-Gaddum Problem is to
determine sharp bounds for (1) $f(G) + f(G^c)$ and (2) $f(G)f(G^c)$,
and characterize the extremal graphs. The Nordhaus-Gaddum type
relations have received wide attention; see a recent survey paper
\cite{Aouchiche-Hansen} by Aouchiche and Hansen. Theorem~\ref{thmf}
concerns such type of a problem for the parameter $\lambda_k$. To
prove the theorem, we will need the following proposition.

\begin{pro}\label{pro8}
A digraph $D$ is strong if and only if $\lambda_k(D)\ge 1$ for
$2\leq k\leq n$.
\end{pro}
\begin{pf} If $D$ is strong, then there is an out-branching rooted at any vertex $r\in V(D)$, so $\lambda_n(D)\ge
1$. By Proposition~\ref{proa}(\ref{pro1}), $\lambda_k(D)\ge 1$ holds
for $k\geq 2$. Now assume that $\lambda_k(D)\ge 1$. By
Theorem~\ref{thmb}, $\lambda(D)\geq 1$, so $D$ is strong.
\end{pf}

Ng proved the following result on the Hamiltonian decomposition of
complete regular multipartite digraphs.

\begin{thm}\label{thm06}\cite{Ng}
The arcs of $\overleftrightarrow{K}_{r, r, \ldots, r}$ ($s$ times)
can be decomposed into Hamiltonian cycles if and only if $(r, s)\not=(4,1)$ and $(r, s)\not=(6, 1)$.
\end{thm}

Now we can get our sharp Nordhaus-Gaddum type bounds for the
parameter $\lambda_k(D)$.
\begin{thm}\label{thmf}
For a digraph $D$ with order $n$, the following assertions hold:
\begin{description}
\item[(i)]~$0\leq \lambda_k(D)+\lambda_k(D^c)\leq n-1$. Moreover, both bounds are sharp. In particular, the lower bound holds if and only if $\lambda(D)=\lambda(D^c)=0$.
\item[(ii)]~$0\leq \lambda_k(D){\lambda_k(D^c)}\leq
\lfloor(\frac{n-1}{2})^2\rfloor$. Moreover, both bounds are sharp.
In particular, the lower bound holds if and only if $\lambda(D)=0$
or $\lambda(D^c)=0$.
\end{description}
\end{thm}
\begin{pf}
\noindent{\bf Part $(i)$}. Since $D\cup
D^c=\overleftrightarrow{K}_n$,
Proposition~\ref{proa}(\ref{pro3}) implies the following, where $x \in V(D)$ is arbitrary.

\[
\lambda_k(D)+\lambda_k(D^c)\leq   \delta^+(D) + \delta^+(D^c) \leq d_D^+(x) + d_{D^c}^+(x) = n-1
\]

If $H\cong \overleftrightarrow{K}_n$, then we have $\lambda_k(H)=n-1$ and
$\lambda_k(H^c)=0$, so the upper bound is sharp. The lower bound is
clear. Furthermore, the lower bound holds, if and only if
$\lambda_k(D)=\lambda_k(D^c)=0$. This is the case if and only if
$\lambda(D)=\lambda(D^c)=0$ by Proposition~\ref{pro8}. As an example,
consider a non-strong tournament, $D$, in which case it is not difficult to see that $D^c$ is also non-strong.

\noindent{\bf Part $(ii)$}. The lower bound is clear. Furthermore,
the lower bound holds, if and only if $\lambda_k(D)=0$ or
$\lambda_k(D^c)=0$ (i.e. if and only if $\lambda(D)=0$ or
$\lambda(D^c)=0$ by Proposition~\ref{pro8}).
For the upper bound, we have

\[
\lambda_k(D){\lambda_k(D^c)}\leq
\left(\frac{\lambda_k(D)+\lambda_k(D^c)}{2}\right)^2\leq
\left(\frac{n-1}{2}\right)^2.
\]

Since both $\lambda_k(D)$ and $\lambda_k(D^c)$ are integers, the upper bound holds.

We now prove the sharpness of the upper bound. By
Theorem~\ref{thm06}, the complete regular bipartite digraphs
$\overleftrightarrow{K}_{a, a}$ with bipartite sets $A=\{u_i\mid
1\leq i\leq a\}$ and $B=\{v_i\mid 1\leq i\leq a\}$ can be decomposed
into $a$ Hamiltonian cycles: $D_i~(1\leq i\leq a)$. We could relabel
the vertices of $\overleftrightarrow{K}_{a, a}$ such that $D_a=u_1,
v_1, \ldots, u_i, v_i, \ldots, u_a, v_a, u_1$.

Let $D$ be the union of the former $a-1$ Hamiltonian cycles: $D_i~(1\leq
i\leq a-1)$. Then there are $a-1$ pairwise arc-disjoint
out-branchings rooted at any vertex of $D$, so $\lambda_n(D)\geq
a-1$. By Proposition~\ref{proa}(\ref{pro1}), $\lambda_k(D)\geq a-1$
for any $2\leq k\leq n$. Furthermore, since
$\delta^+(D)=\delta^-(D)=a-1$, by
Proposition~\ref{proa}(\ref{pro3}), we have $\lambda_k(D)=a-1$ for
any $2\leq k\leq n$. Now $D^c$ is a union of a complete digraph with
vertex set $A$, a complete digraph with vertex set $B$ and the
Hamiltonian cycle $D_a$. We now compute the maximum number of
pairwise arc-disjoint out-branchings rooted at any vertex $r$ in
$D^c$. Without loss of generality, assume that $r=u_1$. By
Lemma~\ref{thm1}, in $D^c[A]$, there are $a-1$ pairwise arc-disjoint
out-branchings rooted at $r$: $T'_i~(1\leq i\leq a-1)$. For $1\leq i\leq a-1$, let
$T_i$ be union of $T'_i$ and the arc set $\{u_{i+1}v_{i+1},
v_{i+1}v\mid v\in B\setminus \{v_{i+1}\}\}$; let $T_a$ be an
out-branching with the arc set $\{u_{1}v_{1}, v_{1}v\mid v\in
B\setminus \{v_{1}\}\}\cup \{v_iu_{i+1}\mid 1\leq i\leq a-1\}$. It
can be checked that $T_1, \ldots, T_a$ are $a$ pair-wise
arc-disjoint out-branchings rooted at $r$, so 
$\lambda_n(D^c)\geq a.$ By Proposition~\ref{proa}(\ref{pro1}),
$\lambda_k(D^c)\geq a$ for any $2\leq k\leq n$. Furthermore, since
$\delta^+(D^c)=\delta^-(D^c)=a$, by
Proposition~\ref{proa}(\ref{pro3}), we have $\lambda_k(D^c)=a$ for
$2\leq k\leq n$. Hence, $\lambda_k(D){\lambda_k(D^c)}=a(a-1)$. Now
in this case the upper bound is
$\lfloor(\frac{2a-1}{2})^2\rfloor=\lfloor(a^2-a+1/4)\rfloor=a^2-a=a(a-1)=\lambda_k(D){\lambda_k(D^c)}$.
This completes the proof.
\end{pf}

\section{Discussions}

In this paper, we completely determine the complexity for both $\kappa_{S, r}(D)$ and $\lambda_{S, r}(D)$ on general digraphs, symmetric digraphs and Eulerian digraphs, as shown in Tables~1-6. We have not considered further classes of digraphs, but it would be interesting to determine the complexity for other classes of digraphs, like semicomplete digraphs. For example, one may consider the following question.

\begin{op}\label{op2}
What is the complexity of deciding whether $\kappa_{S,r}(D)\ge \ell~$ (respectively, $\lambda_{S,r}(D)\ge \ell)$
for integers $k\ge 3$ and $\ell\ge 2$, and a semicomplete digraph $D$?
\end{op}

Recall that in the argument for the question of deciding whether $\kappa_{S,r}(D)\ge \ell~$ for a symmetric digraph $D$
when both $k$ and $\ell$ are fixed, we use Corollary~\ref{corLinkageUndirected} where the $2k$ vertices
$s_1, s_2, \ldots, s_k, t_1, t_2, \ldots, t_k$ are not necessarily distinct.
Normally in the $k$-linkage problem the initial and terminal vertices are considered distinct.  However if we allow them to be non-distinct and look
for internally disjoint paths instead of vertex-disjoint paths, and the problem remains polynomial, then a similar approach to
that of Theorem~\ref{thmPOLsym} and Corollary~\ref{corPOLsym} can be used to show polynomiality of deciding whether $\kappa_{S,r}(D)\ge \ell~$.

\vskip 1cm

\noindent {\bf Acknowledgement.}  We are thankful to Professor Gregory Gutin for discussions on the concept of directed tree
connectivity when Yuefang Sun visited Royal Holloway University of London, and Professor Magnus Wahlstr\"{o}m for
discussions on the complexity of linkage problem for Eulerian digraphs.

\end{document}